\documentclass[a4paper]{article}
\usepackage[affil-it]{authblk}
\usepackage{graphicx}
\usepackage[space]{grffile}
\usepackage{latexsym}
\usepackage{amsfonts,amsmath,amssymb}
\usepackage{url}
\usepackage[utf8]{inputenc}
\usepackage{hyperref}
\hypersetup{colorlinks=false,pdfborder={0 0 0}}
\usepackage{textcomp}
\usepackage{longtable}
\usepackage{multirow,booktabs}

\usepackage{amscd,theorem}
\usepackage{color}
\usepackage[all]{xypic}
\usepackage[english]{babel}
\makeatletter
\newcommand{\specificthanks}[1]{\@fnsymbol{#1}}

\usepackage[all]{xypic}

\newtheorem{theorem}{Theorem}
\newtheorem{teo}{Theorem}

\newtheorem{lema}[theorem]{Lemma}

\newtheorem{defi}{Definition}[section]
\newtheorem{ejem}{Example}[section]

\numberwithin{equation}{subsection}

\def\proof{\medskip \noindent \textit{Proof: }}
\def\qed{\hfill $\square$ }
\def\Ker{\mathrm{Ker}\,}

\def\T{\mathsf{T}}

\def\LL{\mathcal{L}}

\def\d{\mathrm{d}}

\def\R{\mathbb{R}}

\def\EE{\mathcal{E}}

\def\VolX{{\rm Vol}}
\def\RR{\mathbb{R}}

\def\NatTen{\T }
\def\UnivTen{\T }
\def\Tensors{ p\mbox{-}Tensors}

\def\Sk{\mathsf{S}_{2,k}}
\def\K{\mathsf{K}}
\def\GB{\mathsf{P}}

\begin{document}

\title{On the uniqueness of the Gauss-Bonnet-Chern formula \\ (after
Gilkey-Park-Sekigawa)}


\author{Navarro, A.\\ ICMat, Madrid, Spain \and Navarro,
J.\thanks{This author has been partially supported by Junta de Extremadura and FEDER funds}\\ Universidad de
Extremadura\\Badajoz, Spain \\ navarrogarmendia@unex.es}

\date{}

\maketitle 

\begin{abstract}
On an oriented Riemannian manifold, the Gauss-Bonnet-Chern formula asserts that the Pfaffian of the metric represents, in de Rham cohomology, the Euler class of the tangent bundle. Hence, if the underlying manifold is compact, the integral of the Pfaffian is a topological invariant; namely, the Euler characteristic of the manifold. 

In this paper we refine a result originally due to Gilkey (\cite{Gilkey_1975}) that characterizes this formula as the only (non-trivial) integral of a differential invariant that is independent of the underlying metric.
\end{abstract}

\tableofcontents

\section{Introduction}

In what follows, $\,X\,$ will always denote a compact and oriented smooth manifold.

Let $\,g\,$ be a pseudo-Riemannian metric on $\,X$ and let $\,R\,$ denote its curvature tensor. If we consider the former as a 1-form with values on 1-forms, and the latter as a 2-vector with values on 2-vectors, then
$$ g^{2k} := g \wedge \stackrel{2k }{\ldots} \wedge g \quad \mbox{ and } \quad R^k := R \wedge \stackrel{k}{\ldots} \wedge R $$ are sections of $\,\Lambda^{2k} T^* X \otimes \Lambda^{2k} T^*X\,$ and $\,\Lambda^{2k} TX \otimes \Lambda^{2k} TX$ respectively.

If $\,c\,$ stands for the contraction operator:
$$ \left( \Lambda^{2k} TX \otimes \Lambda^{2k  } T X \right) \otimes \left(
\Lambda^{2k} T^* X \otimes \Lambda^{2k} T^* X \right)
\ \xrightarrow{ \quad c \quad }  \
\R _X\  , $$ that matches the first $2k$-indices of the first factor with the first $2k$-indices of the second one, and the last $2k$-indices with last $2k$-indices, then, for any $\,k\geq 0$, we can consider the following functions:
$$ \GB_k := \frac{1}{(2\pi)^k k!}  \, c \left( R^k \otimes g^{2k} \right) \ . $$ 

In even dimenion $\,n=2k$, the top form $\,\GB_k \VolX_g\,$, where $\,\VolX_g\,$ denotes the volume form associated to the metric $\,g$, is called the Pfaffian of the metric.

Then, the so-called Gauss-Bonnet-Chern formula states that the Pfaffian $\, \GB_k \VolX_g\,$ represents, in de Rham cohomology, the Euler class of $\,TX$, for any Riemannian choice of $\,g$ (see \cite{Milnor_1974}). 

When combined with the Poincar\'{e}-Hopf theorem, there is an integral version of this fact, whose statement is more convenient for our purposes (see \cite{Gilkey_2015}):

\medskip
\noindent {\bf Theorem (Gauss-Bonnet-Chern):}
{\it Let $\,X\,$ be a compact and oriented smooth manifold of even dimension $\,n= 2k$.

For any pseudo-Riemannian metric $\,g\,$ on $\,X$ with signature $\,(n_+ , n_{-})$, it holds:
$$ \int_X \GB_k (g) \VolX_g   \, = 
\begin{cases} 
\, 0 & \mbox{ if } n_{-} \mbox{ is odd,} \\
(-1)^{n_{-} /2}\, \chi (X) & \mbox{ if } n_{-} \mbox{ is even,} 
\end{cases} $$  where  $\,\chi (X)\,$ denotes the Euler characteristic of $\,X$.
}
\medskip

The problem of the uniqueness of this formula was already considered in \cite{Gilkey_1975} and can be formulated in the following terms:

Let $\, Metrics\,$ and $\,Functions\,$ denote the sheaves of pseudo-Riemannian metrics and functions over $\,X$, respectively. A morphism of sheaves $\,\LL \colon Metrics \to Functions\,$ is called a {\it differential invariant} if it is smooth and natural (Definition \ref{definicionNatTensor}). It is homogeneous of weight $\, w \in \mathbb{Z}\,$ if $\, \LL (\lambda^2 g) = \lambda^w \LL (g)$, for any metric $\,g$ and any real number $\, \lambda \in \RR$. For example, the functions $\,\GB_k\,$ defined above are differential invariants, homogeneous of weight $-2k$.

For any fixed dimension, the vector space of differential invariants is independent of the underlying manifold; i. e., any such a $\,\LL\,$ is is defined over any smooth manifold of the fixed dimension.

So, let $\,\LL\,$ be a homogeneous differential invariant in dimension $\,n$, and suppose the following hypothesis:
\begin{itemize}
\item[(H)] There exists a compact and oriented smooth manifold $\,X\,$ of dimension $\,n\,$ such that:
$$ \int_{X} \LL (g) \VolX_g  \mbox{ is a non-zero number, independent of the metric } g  \mbox{ on } X . $$
\end{itemize}
\
Then, the question is whether it necessarily follows the existence of $\,\mu \in \RR$ such that:
$$\, [\LL (g)\VolX_g ] \, =\, \mu \, [\GB_k (g) \VolX_g]  $$ as cohomology classes in $\,H^{n} (X, \RR)$ (i. e., modulo exact forms), for any metric $\,g$.

This question was first solved affirmatively by P. Gilkey in \cite{Gilkey_1975}, where he assumed that $\,\LL\,$ was a polynomial differential operator; i. e., a polynomial function $\, \LL \colon J^k M \to \RR\,$ on the space of $k$-jets of metrics, for some $\, k\in \mathbb{N}$. Later on, Gilkey himself relaxed this hypothesis (\cite{Gilkey_1978}), by simply requiring $\,\LL\,$ to be a smooth function $\,\LL \colon J^kM \to \RR$, for some $\, k\in \mathbb{N}$.

In this paper, we go another step further in this direction, by only assuming $\,\
\LL\,$ to be a morphism of sheaves $\, \LL \colon Metrics \to Functions$. 

Thus, we prove:

\medskip
\noindent {\bf Theorem \ref{Principal}:} {\it Let $\,\LL\colon Metrics \to Functions\,$ be a homogeneous differential invariant.

Assume:

\begin{itemize}
\item[(H)] There exists a compact and oriented smooth manifold $\,X\,$ of dimension $\,n\,$ such that:
$$ \int_{X} \LL (g) \VolX_g  \mbox{ is a non-zero number, independent of the metric } g  \mbox{ on } X . $$
\end{itemize}

Then, $\,n = 2k\,$ is even, and there exists $\,\mu \in \R\,$ such that:
$$ [ \LL (g)\VolX_g ] \, =\, \mu \, [ \GB_k (g) \VolX_g ] \ , $$ for any metric $\,g\,$ on any smooth manifold of dimension $\,n$.}

\medskip
To prove this result, we use standard machinery from the theory of natural operations (\cite{Navarro_2008},\cite{Navarro_2015}), as well as recent computations regarding certain dimensional identities of the curvature (\cite{Gilkey_2012}, \cite{Navarro_2014}).

Finally, let us also mention that the proof presented here, that follows ideas initiated by Gilkey-Park-Sekigawa, takes advantadge of the relation between dimensional curvature identities and variationally trivial lagrangians. The long-standing Deser-Schwimmer conjecture, which is a similar problem to that addressed here but where hypothesis $\,(H)\,$ is modified by saying that $\, \int_{X} \LL (g) \VolX_g\,$ is independent of the {\it conformal class} of $\, g$, was finally established with an overwhelmingly lengthy proof (\cite{Alexakis_2009}, \cite{Alexakis_2012}). In our opinion, the ideas of Gilkey-Park-Sekigawa may also shed some light on this topic.

\section{Universal tensors}

Fix a dimension $n \in \mathbb{N}$. 

Let $\,X\,$ be a smooth manifold of dimension $\,n$, and let $\,\Tensors\,$ and $\,Metrics\,$ denote the sheaves of $p$-covariant tensors and pseudo-Riemannian metrics of a fixed signature $\,(n_+ , n_{-})$, respectively.

\begin{defi}\label{definicionNatTensor}
A  $p$-covariant, {\bf  natural  tensor in dimension $n$}, {\bf associated to metrics of signature $\,(n_+ , n_{-})$}, is a morphism of sheaves
$$ T \colon Metrics \longrightarrow \Tensors $$ such that:
\begin{itemize}
\item It is {\it smooth}: if $\{ g_t \}_{t \in T}$ is a smooth family of metrics on an open subset $\,U \subset X$, parametrized by a smooth manifold $T$, then $\{ T(g_t) \}_{t\in T}$ is also a smooth family of tensors.

\item It is {\it natural}: for any local diffeomorphism $ \tau \colon U \to V$ between open sets of $\,X$, it holds:
$$
T(\tau^* g)= \tau^* T (g) \ . $$
\end{itemize}

A natural tensor $ T $ is \textbf{homogeneous of weight} $w \in \mathbb{R}$ if, for any metric $g$ and any positive real number $\lambda>0$, it holds:
$$
T ( \lambda^2 g ) = \lambda^{w}\, T ( g)  \ . $$
\end{defi}

Natural tensors do not depend on the underlying manifold $\,X$, neither on the fixed signature of the metrics: the vector space of homogeneous 
natural tensors associated to Riemannian metrics in $\,\RR^n\,$ is
canonically isomorphic to the space of homogeneous natural tensors associated to non-singular
metrics of signature $\,(n_{+} , n_{-})\,$ over an arbitrary smooth $\,n$-dimensional manifold $\,X$.

Consequently, the $\RR$-vector space of homogeneous natural tensors in dimension $n$ will be denoted:
$$
\NatTen _{p,w} [n] :=   \left[
\begin{array}{c}
\text{ $p$-Covariant, natural tensors $T$ in dimension $n$ } \\
\text{ homogeneous of weight } w \
\end{array}
\right] \ . $$

\begin{ejem}
If $\,p=0$, natural $0$-tensors are also called differential invariants. The expressions $\,\GB_k\,$ defined in the introduction are examples of differential invariants:
$$ \GB_k \colon Metrics \to Functions \ . $$

Observe that each of them is defined for any dimension $\,n$, so that they are in fact natural in a slightly stronger sense, that is explained below. 
\end{ejem}

Let $(X,g)$ be an $(n-1)$-Riemannian manifold, and consider the cylinder $(X\times \RR \, , \, g + \d t^2)$, which is an $n$-dimensional Riemannian manifold.

Let $i$ denote the embedding by the equator:
$$ i \colon X \hookrightarrow  X\times \RR \qquad , \qquad x \ \mapsto \ (x,0) \ . $$

\begin{defi}
The \textbf{dimensional
reduction} of natural tensors is the linear map:
$$
\NatTen _{p,w} [n] \xrightarrow{\ \ r_n \ \ } \NatTen _{p,w} [n-1] \quad , \quad  r_n(T)(g):=i^* \left( T(g + \d t^2) \right) \ .
$$ \end{defi}

If $T \in \NatTen _{p,w} [n] $, it is not difficult to check that $r_n(T)$ is a natural tensor in dimension $n-1$. 

Moreover, as $\lambda^2 g + \d t^2$ and $\lambda^2 g + \lambda^2 \d t^2$ are related by an isometry that is the identity on $X$,
\begin{align*}
r_n(T) (\lambda^2 g) &= i^* \left( T(\lambda^2 g + \d t^2 ) \right) = i^* \left( T ( \lambda^2 g + \lambda^2 \d t^2 ) \right) \\
& = i^* \left( \lambda^w T( g + \d t^2) \right) = \lambda^w \, r_n(T) ( g)  \ ,
\end{align*} so that $\,r_n(T)\,$ is also homogeneous of weight $\,w$.

Therefore, these dimensional reduction maps establish a projective system:
$$
\ldots \xrightarrow{\quad r_{n+1} \quad} \NatTen _{p,w}[n] \xrightarrow{\quad r_n \quad}
\NatTen _{p,w}[n-1]  \xrightarrow{\quad r_{n-1}\quad} \ldots \ .
$$

\begin{defi}
A $p$-covariant {\bf universal tensor}, homogeneous of weight $w$, is an element of the inverse limit $$ \NatTen _{p,w}:= \varprojlim \NatTen _{p,w}[n] \ . $$
\end{defi}

\medskip

In other words, a universal tensor is a collection of natural
tensors $\, \{ T_n \}_{n\in \mathbb{N}} $,
where each $\,
T_{n}\,$ is a natural tensor in dimension $n$, satisfying that, for
every manifold $\,X\,$ of dimension $n_1$ and every embedding into a cylinder
$i \colon X \rightarrow X \times \mathbb{R}^{n_2}$, $x \ \mapsto \ (x,0)$, where $\mathbb{R}^{n_2}$ is endowed with the euclidean metric, it holds
$$i^*\left( T_{n_1 + n_2} ( g + \sum^{n_2}_{i=1} \d t_i^2 ) \right) = T_{n_1}(g) \ . $$

\medskip

\begin{ejem}
For any fixed $\lambda\in \RR$, the tensor $\lambda g$ is  a universal tensor. However,  $ ({\rm tr} \, {\rm Id}) \, g$, $(-1)^{\dim X} g$ or $(-1)^{n_+} g$  are not universal tensors.

The Riemann-Christoffel tensor $R$, the Ricci tensor or the scalar curvature are also examples of universal tensors.
\end{ejem}

\begin{ejem}\label{Pfaffian}
For any $k\geq 0$, the differential invariant defined in the introduction:
$$ \GB_k := c \left( R^k \otimes g^{2k} \right) \  $$ is universal, for it is obtained contracting indices of the curvature with indices of a universal tensor (\cite{Navarro_2014}, Lemma 1.7). 

Its local expression may be written as follows:
$$
(\GB_k) _{i_1 \ldots i_{2}} :=
R^{a_1a_2,b_1b_2} \ldots R^{a_{2k-1}a_{2k}, b_{2k-1}b_{2k}}  g_{a_1 c_1} \ldots g_{a_{2k} c_{2k}} \delta_{b_1 \ldots
b_{2k}}^{c_1 \ldots c_{2k}}
$$
where $\delta_{i_1 \dots i_m}^{j_1 \ldots j_m}$ denotes the generalized Kronecker delta.
\end{ejem}

\begin{ejem}
With the notations of the above example, let us define, for any $k\geq 0$, the $2$-covariant universal tensor:
$$ \Sk := c \left( R^k \otimes g^{2k + 1} \right) \ . $$


Analogously, each $\Sk$ is indeed a universal tensor, with local expression:
$$
(\Sk) _{i_1 i_{2}} :=
R^{a_1a_2,b_1b_2} \ldots R^{a_{2k-1}a_{2k}, b_{2k-1}b_{2k}}  g_{a_1 c_1} \ldots g_{a_{2k} c_{2k}} g_{j i_{2}} 
\delta_{b_1 \ldots b_{2k} i_1}^{c_1 \ldots c_{2k} j} \ . 
$$
\end{ejem}

\subsection{Dimensional identities}

Let us now focus our attention on natural tensors with $\,p=2\,$ indices. 
First of all, let us recall that, if $\,T \in \NatTen _{2  , w}\, $ is such a
2-covariant tensor, homogeneous of degree $\,w$, then $\,w\,$ has to
be an even integer, lesser or equal than 2 (\cite{Navarro_2014}, Section 1.2.1); that is:
$$ w = 2  - 2k $$ for some $k \geq 0$.

On the other hand, a 2-covariant, {\it curvature identity} in dimension $n$, homogeneous of weight $\,w\,$ is an element of the projective space associated to the following vector space:
$$ \K _{2 , w} [n] := \Ker \left[ \phantom{\frac{1}{1}} \hskip -.3cm  \NatTen _{2,w}  \longrightarrow \NatTen _{2,w} [n] \right]   \ . $$

We will make essential use of the recent computation of some of these kernels:

\begin{teo}[\cite{Gilkey_2012},\cite{Navarro_2014}]\label{IdentidadesDimensionales}
For any fixed weight $w = 2 -2k$, with $k \geq 1$, the following holds:

\begin{itemize}
\item If $\,n \geq 2k + 1$, there are no dimensional curvature identities; i.e.:
$$ \K_{2, w} \left[ n \right] \, = \, 0 \quad , \quad \mbox{ for } n \geq  -w + 3 \ .  $$

\item If $\,n = 2k $, the only dimensional curvature identity is the vanishing of $\mathsf{S}_{2, k} $:
$$ \K_{2, 2-n} \left[ n \right] \, = \, \langle\mathsf{S}_{2, k}  \rangle \quad , \quad \mbox{ for } n = -w + 2 \ .  $$
\end{itemize}
\end{teo}

\section{Some lemmas from variational calculus}

Let $\,X\,$ be oriented, and let $ \, \VolX_g \,$ denote the volume form associated to a metric $\,g$.




Any differential invariant $\, \LL \colon Metrics \to Functions \,$ can be plugged into a {\it lagrangian density} $\,\LL (g) \VolX_g\,$.

Standard techniques from the calculus of variations (e. g. \cite{Anderson_1984}) allow to construct an {\it Euler-Lagrange tensor} $\,\mathbb{E} \,$ out of any lagrangian density.
As $\,X\,$ is oriented, the Euler-Lagrange tensor of a lagrangian density can be identified with a 2-covariant tensor $\,\mathcal{E} (\LL )_g\,$ on $\,X$:
$$ \mathbb{E} ( \LL(g) \VolX_g ) \, = \, \EE (\LL)_g \otimes \VolX_g \ . $$

The construction of the Euler-Lagrange tensor is natural and homogeneous, in the sense that the following squares commute (\cite{Anderson_1984}, Lemma 2.3):

\begin{equation*}\label{Extension}
\xymatrix{
\UnivTen_{0,w} [n] \, \ar[r]^-{\mathcal{E}} \ar[d]_-{r_n} & \UnivTen_{2,w+2} [n]  \ar[d]^-{r_n} \\
\UnivTen_{0,w} [n-1] \, \ar[r]^-{\mathcal{E}}   & \, \UnivTen_{2,w+2} [n-1]  }  
\end{equation*} 

That is to say,

\begin{lema}\label{LemaPesoEL}
The Euler-Lagrange construction $\, \mathcal{L} \longmapsto \mathcal{E}(\mathcal{L})\,$ defines a linear map between the space of differential invariants, homogeneous of weight $\,w$, and that of universal $2$-tensors, homogeneous of weight $\,w + 2$:
$$ \mathcal{E} \colon \UnivTen _{0,w} \longrightarrow \UnivTen_{2,w+2}  \ . $$
\end{lema}

Loosely speaking, the value of the Euler-Lagrange tensor $\,\EE (\LL)\,$  can be understood as the differential of the function:

$$
\{ \mbox{ Metrics on } X \, \} \ \longrightarrow \ \RR 
\quad , \quad g \ \longmapsto \ \int_X \LL (g) \VolX _g \ . 
$$

In fact, it is not difficult to prove:

\begin{lema}\label{DiferencialEsCero}
Let $\,\LL \colon Metrics \to Functions\, $ be a differential invariant.

If $\,\int_X \LL (g) \VolX_g \, $ is independent of the metric $\,g\,$ on which the integrand is evaluated, then its Euler-Lagrange tensor is identically zero; that is:
$$ \mathcal{E} (\LL )_g = 0 \quad \forall \ g \  . $$

On the other hand, if there exists a metric such that $\,\int_X \LL (g) \VolX_g  \neq \int_X \LL (\lambda g) \VolX_{\lambda g} \,$ for any sufficiently small $\lambda \in \R^+ $, then $\, \mathcal{E} (\LL)_g \neq 0$.
\end{lema}

Next, we will use a standard lemma from the inverse theory of the calculus of variations (\cite{Anderson_1992}, Thm. 2.7). To state it, let $\,Vectors\,$ denote the sheaf of smooth vector fields on $\,X$.

\begin{lema}\label{MismoEL}
If $\,\LL \,$ and $\,\LL '\,$ are differential invariants with the same Euler-Lagrange tensor, then there exists a smooth morphism of sheaves $\,D\colon Metrics \to Vectors\,$ such that, for any metric $g$:
$$ \LL (g) \ = \ \LL ' (g) \, + \, {\rm div} \left( D(g)\right) \ . $$
\end{lema}

Finally, the following is a classical computation that can be found, for example, in \cite{Lovelock_1971}:

\begin{lema}\label{ELdeLovelock}
For any $\,k \geq 0$, the Euler-Lagrange tensor of $\,\GB_k\,$ is (proportional to) the tensor $\,\Sk$:
$$ \mathcal{E} ( \GB_k )  \, = \, \frac{1}{2} \, \Sk \ . $$
\end{lema}

\section{Characterization of the Gauss-Bonnet-Chern formula}

Let us prove in this Section the announced characterizacion of the Gauss-Bonnet-Chern formula.

Firs of all, let us recall that elements of $\,\UnivTen _{0,w}\,$, are smooth and natural morphisms of sheaves $\, \LL \colon Metrics \to Functions$, homogeneous of weight $\,w$, and which are defined for any dimension of the underlying manifold.

\begin{teo}\label{Principal}
Let $\,\LL \colon Metrics \to Functions \,$ be a homogeneous differential invariant; i. e., an element of $\, \UnivTen _{0,w}$.

Assume:

\begin{itemize}
\item[(H)] There exists a compact and oriented smooth manifold $\,X\,$ of dimension $\,n\,$ such that:
$$ \int_{X} \LL (g) \VolX_g  \mbox{ is a non-zero number, independent of the metric } g  \mbox{ on } X . $$
\end{itemize}

Then, $n = 2k$ is even, and there exist $\,\mu \in \R\,$ and a smooth morphism of sheaves $\,D\colon Metrics \to Vectors \,$ such that:
$$ \LL (g)\, =\, \mu \, \GB_k (g) + {\rm div} \left( D (g) \right) \ , $$ for any metric $\,g\,$ on any smooth manifold of dimension $\,n$.
\end{teo}

\proof First of all, observe that, for any $\lambda \in \R$,
$$ \int_X \mathcal{L}(\lambda^2 g ) \VolX_{\lambda^2 g} = \lambda^{w + n} \int_X \mathcal{L} (g) \VolX_g \ , $$ so that, for the hypothesis to be held, the weight of $\LL$ has to be $\,-n$.

Hence, its Euler-Lagrange tensor $\, \EE (\LL)\,$ is a universal 2-tensor homogeneous of weight $2 - n$ (Lemma \ref{LemaPesoEL}): 
$$ \EE (\LL) \in \T_{2, 2 - n} \ . $$

As the integral $\,\int_X \mathcal{L}(g) \VolX_g \,$ does not depend on the metric, Lemma \ref{DiferencialEsCero} implies that  $\,\EE (\LL)\,$  vanishes on dimension $n$; that is:
$$ \EE (\LL) = 0 \in \T_{2,2-n} [n] \ . $$

On the other hand, if $i_\theta \colon X \to X \times S_1$ denotes the inclusion $x \mapsto (x, \theta)$:
\begin{align*}
\int_{X \times S_1} \LL (g + \d \theta^2) \VolX_{g + \d \theta^2}  &= \int_{S_1} \left( \int_X i^*_\theta ( \LL (g + \d \theta^2 ) \VolX_{g + \d \theta^2 })  \right) \\
& = \int_{S_1} \left( \int_X ( \LL (g ) \VolX_{g  })  \right) \d \theta \\
 &=  2 \pi   \int_X ( \LL (g ) \VolX_{g  }) \, \neq \, 0  . 
\end{align*}

Therefore, when $h = g + \d \theta^2$, 
$$ \int_{X \times S_1}  \LL (\lambda^2 h ) \VolX_{\lambda^2 h} = \lambda^{-n + n + 1}   \int_{X \times S_1}  \LL (h ) \VolX_{h}  \neq \int_{X \times S_1} \LL (h) \VolX_{h} \ , $$
so that,
\begin{equation*}\label{ArgumentoFinal}
\EE (\LL) \neq 0 \ \mbox{ in } \ \T_{2, 2-n}[n+1] \ . 
\end{equation*} 

In other words, $\,\EE (\LL)\,$ belongs to the kernel $\, \K_{2,2-n}[n]$ of the dimensional restriction $\,r_{n+1} \colon \UnivTen_{2,2-n} [n+1] \to \UnivTen_{2,2-n} [n]$. 

And Theorem \ref{IdentidadesDimensionales} precisely computes this kernel: 
$$ \K_{2, 2 - n} \left[ n \right] \, = \, \langle\mathsf{S}_{2, k}  \rangle \ . $$

That is to say, there exists $\,\widetilde{\mu} \in \R\,$ such that the following equality holds in $\,\UnivTen _{2,2-n} [n]$:
$$\, \EE (\LL) \, = \, \widetilde{\mu} \, \Sk \, = \, 2 \widetilde{\mu} \, \EE ( \GB_k ) \ . $$

Finally, Lemma \ref{MismoEL} then assures the existence of a smooth morphism of sheaves $\,D\colon Metrics \to Vectors\,$ such that, for $\,\mu = 2 \widetilde{\mu}$:
$$ \LL \, =\, \mu \, \GB_k + {\rm div} D \ . $$

\qed


\begin{thebibliography}{00}

\bibitem{Alexakis_2009} Alexakis, S.: \emph{On the decomposition of global conformal invariants. I} Ann. of Math. {\bf 170} (2009) 1241–-1306.

\bibitem{Alexakis_2012} Alexakis, S.: \emph{The decomposition of global conformal invariants: some technical proofs II}, Pacific J. Math. {\bf 260} (2012) 1--88. 

\bibitem{Anderson_1984} Anderson, I. M.: \emph{Natural variational principles on Riemannian manifolds}, Ann.
of Math. 120 (1984) 329--370.

\bibitem{Anderson_1992} Anderson, I. M.: \emph{Introduction to the Variational Bicomplex}, in Mathematical Aspects of Classical Field Theory (ed. by M. Gotay, J.Marsden, V.Moncrief), Comptemporary Mathematics Vol. 132, (1992)

\bibitem{Gilkey_1975} Gilkey, P.: \emph{Local invariants of an embedded Riemannian manifold}, Ann. of Math. {\bf 102} (1975) 187--203.

\bibitem{Gilkey_1978} Gilkey, P.: \emph{Smooth invariants of a Riemannian manifold}, Adv. of Math. {\bf 28} (1978) 1--10.

\bibitem{Gilkey_2015} Gilkey, P.; Park, J. H.: \emph{Analytic continuation, the Chern–Gauss–Bonnet theorem,
and the Euler–Lagrange equations in Lovelock theory for indefinite signature metrics}, J. Geom. Phys. {\bf 88} (2015) 88-93


\bibitem{Gilkey_2011}
Gilkey, P.; Park, J.H.; Sekigawa,  K.: \emph{Universal curvature identities}, Diff.
Geom. App., {\bf 62}  (2011) 814-825.

\bibitem{Gilkey_2012}
Gilkey, P.; Park, J.H.; Sekigawa,  K.: \emph{Universal curvature identities II}, J.
Geom. Phys., {\bf 62}  (2012) 814-825.



\bibitem{Lovelock_1971}
Lovelock, D.: \emph{The Einstein tensor and its generalizations}, J. Math. Phys. \textbf{12} (1971) 498--501.

\bibitem{Milnor_1974}
Milnor, J.; Stasheff, J.: \emph{Characteristic classes}, Princeton University Press, New Jersey (1974).


\bibitem{Navarro_2014} Navarro, A.; Navarro, J.: \emph{Dimensional curvature identities on pseudo-Riemannian geometry}, J. Geom. Phys. {\bf 86} (2014) 554--563.

\bibitem{Navarro_2008} Navarro, J.; Sancho, J. B.: \emph{On the naturalness of Einstein's equation}, J. Geom. Phys. {\bf 58} (2008) 1007--1014.

\bibitem{Navarro_2015} Navarro, J.; Sancho, J. B.: \emph{Natural operations on differential forms}, J. Geom. Phys. {\bf 38} (2015) 159--174.

\end{thebibliography}
\end{document}